\documentclass[11pt]{article}
\usepackage{amsmath,amssymb}
\usepackage{graphicx}
\usepackage{hyperref}
\usepackage{authblk}

\newtheorem{thm}{Theorem}
\newtheorem{prp}{Proposition}
\newtheorem{lem}{Lemma}
\newtheorem{cor}{Corollary}
\newtheorem{rem}{Remark}

\newtheorem{ex}{Example}
\usepackage{breqn}

\title{Invariants and areas of Steiner $4$-chains}
\author[1]{G.Bibileishvili}
\author[2]{A.Diakvnishvili}
\affil[1]{Faculty of Business,Technology and Education, Ilia State University, Tbilisi, Georgia}
\date{\today}

\begin{document}
	
	\maketitle
	
	\begin{abstract}
		We are concerned with the invariants of Steiner chains consisting of four circles. In particular, we compute the invariant
		moments of curvatures in Steiner $4$-chains and give two applications of the obtained formulas. Specifically,
		we present an algorithmic feasibility criterion for Steiner $4$-chains and identify the poristic Steiner chains having 
		extremal areas, which yields a generalisation of the main results of a recent paper by K.Kiradjiev.
		The proofs are based on the invariants of Steiner chains described by R.Schwarz and S.Tabachnikov and on the relations
		between the radii of neighbouring poristic circles established by P.Yiu.
	\end{abstract}

    \vspace*{0.25cm}\noindent{\small {\bf Keywords and phrases}:
	{Steiner chain, Soddy circles, Steiner porism, poristic Steiner chains, Descartes circle theorem,
		moments of signed curvatures, area}
	
	\vspace*{0.25cm}\noindent{\small {\bf MSC 2010:} {52C35,  32S40.}}
	
	\section*{Introduction}
	
Motivated by \cite{kir} and relying on \cite{scta}, \cite{yiu}, we compute the invariant moments of curvatures and the extremal values of areas for Steiner chains consisting of four circles. \\

In the first section we compute the invariant moments of curvatures (bends) of Steiner $4$-chains in terms of the radii of parent (Soddy) circles. In the rest of the paper we apply these formulas to solving two natural problems studied in \cite{kir} in some special cases of Steiner $4$-chains. Namely, in the second
section we give an algorithmic feasibility (realizability) criterion for the radii of Steiner $4$-chains and in the third section we determine the $4$-chains having an extremal sum of areas. The last section contains several remarks
on possible generalizations and further research perspectives. \\

Recall that the term {\it Steiner $n$-chain} refers to a sequence of circles $\delta_i, i=1, \ldots, n,$ in Euclidean plane such that the circles with adjacent indices $mod n$ are externally tangent, each of $\delta_i$ is externally tangent to 
a fixed circle $\gamma$ and internally tangent to another fixed circle $\Gamma$ containing $\gamma$ in its interior disk. The nested circles $\gamma$ and $\Gamma$ are called the {\it parent circles} (or {\it Soddy circles}) of the Steiner chain 
considered \cite{ber}. Steiner chains is a classical topic discussed in many papers, in particular, in the context of {\it Steiner porism} (see, e.g., \cite{ped}, \cite{scta}) which is in the focus of our discussion. \\

It should be mentioned that triples of pairwise externally tangent circles are often called "kissing circles". An ancient result of Apollonius states that, for any triple of "kissing circles", there exist exactly two circles which are tangent
to each of the three given circles. These two circles are called the {\it Soddy circles} of kissing circles. The famous Descartes' circle theorem enables one to compute the curvatures of Soddy circles in terms of the curvatures of the three given "kissing circles" \cite{ber}. \\

We are concerned with several problems arising in the context of the aforementioned Steiner porism.
Its special case we are interested in, states that there exists a one-dimensional family of Steiner $n$-chains of circles having the same Soddy circles \cite{scta}. By a way of analogy with the Poncelet porism \cite{ber} the collection of Steiner $3$-chains with fixed Soddy circles will be called {\it poristic Steiner $n$-chains}.
This point of view suggests, in particular, several extremal problems for poristic Steiner $n$-chains with the fixed Soddy circles. A natural extremal problem concerned with the the sum of areas of poristic Steiner $4$-chains was in one special case considered in a prize-winning paper of K.Kiradjiev \cite{kir}. The answer given in \cite{kir} is correct but the reasoning is definitely incomplete. As an application of our approach we give a rigorous proof of the mentioned result from \cite{kir}. Actually, the main aim of the present paper is to give solutions to two
natural problems considered in \cite{kir}. \\

Certain extremal problems for constrained chains of tangent circles have also been discussed in \cite{gaz}, \cite{guvo} but the setting and context in those papers are essentially different from ours. In the last section, we briefly discuss some
further connections and feasible generalizations of our results and approach. \\

			\section*{Invariants of poristic Steiner $4$-chains} \label{invariantmoments}
 Let $n\geq 3$ be a natural number and $(R, r, d)$ be a triple of positive numbers such that $R>r$ and $d^2 = (R - r)^2 - 4qRr,$ where $q =                         \tan^2\frac{\pi}{n}.$ It is known that such a triple defines a {\it poristic system of Steiner $n$-chains} with the parent (Soddy) circles                         having radii $R, r$ and distance between their centres equal to $d$ (see, e.g., \cite{ped}). Denote by $a = \frac{1}{r}, A = -\frac{1}{R}$ the                         so-called {\it signed curvatures} (bends) of the parent circles. As was shown in \cite{scta}, the first $n-1$ moments $I_k^{(n)}$ of curvatures                         $b_j$ of Steiner $n$-chains are invariant in such a poristic family. Hence they can be expressed through radii $R, r$ of the parent Soddy                         circles or, equivalently, through their signed curvatures $a, A$.
 
 For $n=3$, explicit formulas for the first two moments $I_1, I_2$ of bends follow directly from the Descartes circle theorem
 (see, e.g., \cite{scta}):
 \begin{equation} \label{moment31}
 	I_1^{(3)} = \frac{R-r}{2Rr} = \frac{A+a}{2},
 \end{equation}
 \begin{equation} \label{moment32}
 	I_2^{(3)} = \frac{R^2 - 6Rr + r^2}{8R^2r^2} = \frac{A^2 + 6Aa + a^2}{8}.
 \end{equation}
 
 \begin{rem}
 	It should be noted that the formula for $I_2^{(3)}$ given in \cite{scta} is not correct
 	and should be replaced by (\ref{moment32bis}). A detailed proof of (\ref{moment32bis}) is given below.
 \end{rem}

 We are going to show that, for $n=4$, analogous explicit formulas for the first three moments $I_j^{(4)}$ can be obtained using some                         results of \cite{yiu}. It should be noted that this task was not considered neither in \cite{scta} nor in \cite{yiu}. We use a natural                         straightforward approach. Since the first three moments of bends are invariant they can be computed at any poristic chain. We choose a $4$-                             chain for which this computation is especially simple.

 To this end we use two general results on Steiner $n$-chains given below. 
 
 \begin{prp} \label{poristicrange1}
 	For a given pair of Soddy circles with gauge $(R, r, d)$, the minimal and maximal possible values of poristic
 	radii $r_i$ are
 	$$r_* = \frac{R - d - r}{2}, r^* = \frac{R + d - r}{2},$$
 	while the minimal and maximal values of poristic curvatures are
 	$$b_* = \frac{2}{R+d-r}, b^* = \frac{2}{R-d-r}$$
 	respectively. For any $r \in [r_*, r^*],$ the poristic family contains a circle of radius $r$. 
 \end{prp}
 
 Let $l$ be the line through the centres of Soddy circles called the {\it axis of porism}. It is geometrically obvious and
 easy to prove rigorously that the extremal values of poristic radii are attained at the two poristic circles having
 their centres on line $l$. The radii of those circles are readily expressed in terms of the gauge, which yields
 the formulas given in the lemma. The second statement of the above proposition is a direct consequence of the Steiner porism.
 
 It is also easy to see that among poristic Steiner $n$-chains there always are a chain $C^*$ containing a circle with the biggest possible
 radius $r^*$ and another one, $C_*,$ containing a circle with the smallest possible radius $r_*$ given above. Both these chains are called axial 
 $n$-chains and are axi-symmetric, i.e. symmetric with respect to the line through the centres of Soddy circles called the {\it axis of porism}. 
 For even $n,$ the two chains $C_*$ and $C^*$ coincide and yield the so-called {\it axial $n$-chain}. For even $n$, there are also two axi-symmetric 
 chains having only two common points with the axis at which they are tangent to the axis. The latter two symmetric chains will be called the {\it lateral chains} 
 of poristic family.
 
 Notice that, for $n=4,$ we already know the radii and curvatures of the two extremal circles in $C_*^*$. The remaining two circles are also symmetric with 
 the respect to the axis of porism, and their curvatures can be computed using the following general result which is an immediate corollary of Proposition 5 in \cite{yiu}.
 
 \begin{prp} \label{eqYiu} (\cite{yiu})
 	Let $C(u)$ be a poristic circle in a Steiner $n$-chain having radius $u$ and let $v = \frac{1}{u}$ denote its curvature. Then the curvatures 
 	$v_+, v_-$ of the two neighbours of $C(u)$ are the roots of quadratic trinomial
 	\begin{equation} \label{yiuquadratics}
 		\alpha x^2 + \beta x + \gamma,
 	\end{equation}
 	where
\begin{equation} \label{yiucoefficients}
	\begin{split}
	\alpha = (q+1)^2R^2r^2u^2, \\
	\beta = 2(q+1)Rru[(q-1)Rr - (R-r)s], \\
	\gamma = [(q+1)Rr - (R-r)s]^2 + 4Rru^2. \\
	\end{split}
\end{equation}
 \end{prp}
 
 \begin{cor}
 	The sum $v_+ + v_-$ of two neighboring curvatures of $C(u)$ is equal to $-\frac{\beta}{\alpha}.$
 \end{cor}

We now use this result to obtain formulas for the invariant moments $I_k^{(n)}$ for $n=3, 4.$ First, we put $n=3$ and give another
proof of equations (\ref{moment31}), (\ref{moment32}). Notice that in this case $q=3$, and so the coefficients \label{yiucoefficients}
of equation (\ref{yiuquadratics}) take the form
\begin{equation} \label{yiu3}
	\alpha = 16R^2r^2u^2, \beta = 8Rru[2Rr - (R-r)u],
	\gamma = [4Rr - (R-r)u]^2 + 4Rru^2.
\end{equation}
For the chain $C^*,$ we have $u = r^* = \frac{R-r+d}{2}, v = \frac{2}{R-r+d}$ and an easy computation shows
that in this case
$$-\frac{\beta}{\alpha} = \frac{R-r}{Rr}.$$
Adding $v$ to the above sum of two neighbouring curvatures we finally get
\begin{equation} \label{moment31bis}
	I_1^{(3)} = \frac{A+a}{2} = \frac{R-r}{Rr},
\end{equation}
which of course coincides with (\ref{moment31}). Obviously each of the intermediate bends is equal to $\frac{R-r}{2Rr}.$

Summing now the squares of bends $b_j$ we get
\begin{equation} \label{moment32bis}
	I_2^{(3)} = b_1^2 + b_2^2 + b_3^2 = \frac{A^2 + a^2 + 6aA}{8}= \frac{R^2 - 6Rr + r^2}{8R^2r^2},
\end{equation}
which coincides with (\ref{moment32})

\begin{rem}
	It is easy to verify that the same formulas are obtained by computing these two moments
	for the chain $C_*$, which can be considered as a feasible check of the above formulas.
\end{rem}

We take now $n=4,$ in which case $q=1$. Then for the chain $C_*^*$ we have again $u = r^* = \frac{R-r+d}{2},
v = \frac{2}{R-r+d}$ and the coefficients of equation (\ref{yiuquadratics}) take the form:
\begin{equation} \label{yiu4}
	\alpha = 4R^2r^2u^2, \beta = -4Rru^2(R-r),
	\gamma = [2Rr - (R-r)u]^2 + 4Rru^2.
\end{equation}

Inserting $u=\frac{R-r+d}{2}$ it is easy to verify that the curvatures of two intermediate circles are both equal to $\frac{R-r}{2Rr}.$
So we now know all axial curvatures and present them in a lemma.

\begin{lem} \label{axialbends}
	The quadruple of curvatures of circles in $C_*^*$ is
	\begin{equation} \label{axialcurvatures}
		\left(\frac{2}{R-r-d}, \frac{R-r}{2Rr}, \frac{2}{R-r+d}, \frac{R-r}{2Rr}\right) .
	\end{equation}
\end{lem}

Summing the extremal bends, clearing the denominators and using that $d^2 = R^2 - 6Rr + r^2$ we get
$$\frac{2}{R-r-d} + \frac{2}{R-r+d} = \frac{4(R-r)}{(R-r)^2 - d^2} = \frac{R-r}{Rr}.$$
Adding two intermediate curvatures we finally get the formula for the first moment
\begin{equation} \label{moment41}
	I_1^{(4)} = \frac{2(R-r)}{Rr} = 2(A+a).
\end{equation}

Summing the squares of the above curvatures in a similar way we get the formula for the second moment
\begin{equation} \label{moment42}
	I_2^{(4)} = \frac{3R^2 - 10Rr + 3r^2}{2R^2r^2} = \frac{3A^2 + 10Aa + 3a^2}{2}.
\end{equation}

Finally, summing the cubes of bends we get the following formula for the third moment:
\begin{equation} \label{moment43}
	I_3^{(4)} = \frac{5R^3 - 27R^2r + 27Rr^2 - 5r^3}{4R^3r^3} = \frac{5A^3 + 27A^2a + 27Aa^2 + 5a^3}{4}.
\end{equation}

These explicit formulas for the first three invariant moments of Steiner $4$-chains
will be used below to solve the two problems mentioned in the Introduction.

They also enable us to find the curvatures and radii in the lateral symmetric configuration.
Denote the two unknown lateral curvatures by $b_+$ and $b_-$. From the invariance of the
first two moments of curvatures and the above formulas for $I_1^{(4)}, I_2^{(4)}$ we get
the following system of equations for $b_+, b_-:$
$$\{b_+ + b_- = A+a, b_+^2 + b_-^2 = \frac{3A^2 + 10Aa + 3a^2}{4}\}.$$
It follows that $b_+, b_-$ are the roots of quadratic equation
$$x^2 - (A+a)x + \frac{(A-a)^2}{8} = 0.$$
Hence 
$$b_{\pm} = \frac{A+a}{2} \pm \frac{1}{2\sqrt{2}} \sqrt{A^2 + 6Aa + a^2}.$$
Taking into account that, by Pedoe formula for $n=4,$ we have
$$A^2 + 6Aa + a^2 = \frac{d}{Rr}$$
we finally get 
$$b_{\pm} = \frac{R-r}{2Rr} \pm \frac{d}{2Rr\sqrt{2}}.$$

For convenience of reference in the sequel, we formulate this result as another lemma.

\begin{lem} \label{lateralbends}
	The quadruple of curvatures of circles in both lateral symmetric configurations is
	\begin{equation} \label{offaxialcurvatures}
		\left(\frac{R-r}{2Rr} - \frac{d}{2Rr\sqrt{2}}, \frac{R-r}{2Rr} + \frac{d}{2Rr\sqrt{2}},\\ 
		\frac{R-r}{2Rr} + \frac{d}{2Rr\sqrt{2}}, \frac{R-r}{2Rr} - \frac{d}{2Rr\sqrt{2}}\right).					
	\end{equation}
\end{lem}

In the rest of the paper we apply the above formulas to solving the two problems
mentioned in the Introduction. \\

\section*{\textbf{Feasibility criterion for the radii of Steiner $4$-chains}} \label{feasibility}

We consider first the following feasibility (realizability) problem.
Given an ordered quadruple of positive numbers $(r_1, r_2, r_3, r_4)$ find out if there exists 
a Steiner $4$-chain with such radii in the given order. This natural problem was considered, in particular, 
in \cite{kir}, where the author gave its solution in a concrete case but did not consider the general case.
We give an algorithmic solution to this problem using a different method based on
the formulas obtained above.

Our solution involves four steps. First, we compute the first three moments $(I_1, I_2, I_3)$ of the given
quadruple $r_i$ and call them the {\it actual moments of curvatures}. Then assuming the existence of a sought 
Steiner quadruple we use the equations (\ref{moment41}) and (\ref{moment42}) to find the "virtual radii" 
$(\tilde R, \tilde r)$ of Soddy circles. Thus we wish to solve the system
$${2(a+A)= I_1, 3A^2 + 10Aa + 3a^2 = 2I_2},$$
which reduces to the following quadratic equation
\begin{equation} \label{Soddyradii}
	16a^2 - 8I_1a + (8I_2 - 3I_1^2) = 0.
\end{equation}
which can be solved explicitly. In this way we find the virtual Soddy radii. If $\tilde R^2 - 6 \tilde R \tilde r + \tilde r^2 > 0$ 
then the pair $(\tilde R, \tilde r)$ is a feasible candidate for the radii of Soddy circles with the distance between their centres
equal to $\tilde d = \sqrt{\tilde R^2 - 6 \tilde R \tilde r + \tilde r^2 }.$ Next, we verify that the given values 
of radii belong to the segment $[\tilde r_*, \tilde r^*]$ and use $(\tilde R, \tilde r)$ to compute the first three
"virtual moments of curvatures" $\tilde I_k$ by our formulas (\ref{moment41}), (\ref{moment42}), (\ref{moment43}).
Another necessary condition of realizability is that $I_3 = \tilde I_3).$ Finally, we take any of the given radii, say
$r_1,$ and compute the radii $r_{\pm}$ of its neighbours in the porism with the gauge $(\tilde R, \tilde r, \tilde d$). 
The final condition of our feasibility criterion is that the radii $r_{\pm}$ coincide with $r_4, r_2.$ In such a case we
say that the given quadruple of radii satisfies the {\it feasibility test}. Now we can formulate the criterion we aimed at.

\begin{thm} \label{feasibility4}
	The sought Steiner $4$-chain with radii $(r_1, r_2, r_3, r_4)$ in the given order exists if and only if the quadruple
	$(r_1, r_2, r_3, r_4)$ satisfies the feasibility test.
\end{thm}

{\bf Proof.} The necessity of our feasibility criterion follows directly from the properties of Steiner chains established above.
To prove sufficiency let us take a pair of embedded circles with the gauge $(\tilde R, \tilde r, \tilde d)$. By Pedoe criterion this pair
supports poristic Steiner $4$-chains. Since the given radii belong to the segment $[\tilde r_*, \tilde r^*]$, by Proposition \ref{poristicrange1} 
this poristic family contains circles with each of the given radii. It remains to verify that they come in the given order. This follows
from the fact that the neighboring radii of one them are the same that its neighbors in the given quadruple. The remaining radius will fit
the right place due to the equality between $I_3$ and $\tilde I_3$, which completes the proof.

{\bf Example 1.} For a quadruple of radii $(3, 2.4, 2, 2.4)$, the curvatures are \\ $(0.3333, 0.4166, 0.5, 04166)$ and we get $(I_1, I_2, I_3) = (1.6662, 0.7079, 0.3065).$ The system for virtual Soddy radii
$${2(a+A)= I_1, 3A^2 + 10Aa + 3a^2 = 2I_2},$$
takes the form
$${2(a+A)= 1.662, 3A^2 + 10Aa + 3a^2 = 1.4156}.$$
So in this case equation (\ref{Soddyradii})
$$16a^2 - 8I_1a + (8I_2 - 3I_1^2) = 0$$
takes the form
$$16a^2 - 13.29a - 2.71 = 0.$$

Solving it and taking into account that $a$ should be positive we get that $a=1, A = -1/6.$ So the virtual Soddy radii are $(6, 1)$. 
For these values of Soddy radii, one has $[\tilde r_*, \tilde r^*] = [2,3]$ which contains the given radii. We now compute the adjacent 
radii of the poristic circle with radius $r_1 = 3$ using equation (\ref{eqYiu}) and get that they are equal to $2.4$ and $2.4$, i.e.,
they coincide with the adjacent given radii. Finally, we compute the third virtual moment of curvatures $\tilde I_3$ using 
(\ref{moment43}) and get $\tilde I_3 = 0.3065$, which coincides with $I_3.$ So by our Theorem \ref{feasibility4} this quadruple 
of radii is realizable in a Steiner $4$-chain, which can also be verified geometrically, say using Geogebra.

{\bf Example 2.}  For a quadruple of radii $\rho = (1, 2, 3, 4)$, one can use the same reasoning to show that our criterion is not
fulfilled so these radii are not realizable by any Steiner $4$-chain. Here we have $\varkappa = (1, 0.5, 0.33(3), 0.25)$
and $(I_1, I_2, I_3) = (2.083, 1.4214, 1.1766).$ 	
So here equation (\ref{Soddyradii})  
$$16a^2 - 8I_1a + (8I_2 - 3I_1^2) = 0$$
takes the form
$$16a^2 - 16.64a + 11.37 - 13.02 = 0.$$
Solving it we find out that $\tilde I_3 = 2.305$ is not equal to $I_3 = 1.1766$ so the quadruple $(1, 2, 3, 4)$ is non-realizable.

\begin{rem}
	Actually an analogous algorithm seems feasible for arbitrary $n$. To this end we just need general formulas for the first two moments of poristic curvatures. For even $n$, such formulas can be obtained by computing the curvatures in the axial $n$-chain using iteration of Proposition \ref{eqYiu} starting with the known extremal curvatures of $C_*^*$. 
	For $n=6,$ we give more details in the last section. 
\end{rem} 

For $n=4$, we can get a bit more general result using the general theory of circle packings (see, e.g. \cite{ste}). To this end
recall that if we have a chain of externally tangent circles $\delta_i$ with radii $r_i$ such that there exists a circle $\gamma$
externally tangent to all of them then its radius $r_*$ is algebraically expressible through the radii $r_i$ \cite{ste}. Moreover,
for $n=4$, radius $r_*$ can be found as the inverse of a root of a reduced quartic equation

\begin{dmath}
16x^4-8(b_1b_2 + b_2b_3 + b_3b_4 + b_4b_1 + 2b_1b_3 + 2 b_2b_4)x^2 - 16(b_1b_2b_3 + b_2b_3b_4 + b_3b_4b_1 + b_4b_1b_2)x - 12b_1b_2b_3b_4 - 2(b_1b_2+b_3b_4)(b_2b_3 + b_4b_1) + (b_1^2+b_3^2)(b_2^2 + b_4^2)=0,
\end{dmath}
where $b_i = 1/r_i$ are the curvatures of circles  $\delta_i.$

This gives a necessary condition for the existence of the circle $\gamma$ which is sometimes called the {\it socle}. An analogous equation can 
be written for the radius of circle $\Gamma$ enclosing and internally tangent to all $\delta_i$. So for the virtual radii of Steiner chains we 
get two equations giving the necessary conditions of feasibility. 

A related problem is to find sufficient conditions for the existence of the "socle" $\gamma$ for a given sequence of cyclically
externally tangent circles $\delta_i$. This problem is reduced to solving a system of quadratic equations for the center and radius
of the sought socle circle. 

These problems have connections with the theory of polygonal linkages. Notice that the centers of circles $\delta_i$ form a polygon $P$ with the 
sidelengths equal to the sums of two consequent radii which suggests that it can be considered as a planar polygonal linkage $L(P)$.  
Let us consider in some detail the case where $n=4$. Then polygon $P$ is a tangential quadrilateral and its configuration space as a linkage $L(P)$
is one-dimensional. It is also well known that each convex configuration of $L(P)$ has an incircle. The given circles are in fact its {\it Soddy circles,}
i.e. a system of cyclically tangent circles with the centers at vertices of $P$. Denoting its sides by $(a_1, a_2, a_3, a_4)$ we have 
$a_1 = r_1 + r_2,$ and so on. Let us call this system of circles $P$-necklace and introduce the {\it configuration space of necklace} as all configurations 
of linkage $L(P)$ for which the circles of necklace do not intersect. One can visualize this space as all possible positions of four 
cyclically tangent hard disks (coins), which shows that generically this space is homeomorhic to a closed interval. Steiner feasibility means that
among these configurations there is (at least) one having a socle. A necessary condition is that the above quartic equation has a
positive root on the configuration space of necklace, which in concrete examples can be done using the Sturm algorithm.  

\section*{\textbf{Extremal problems for poristic Steiner $4$-chains}}

Let $F$ be a pair of Soddy circles with the gauge $(R, r, d),$ supporting poristic Steiner $4$-chains and let $b_i$
denote the curvatures of a certain poristic Steiner $4$-chain of circles $\Gamma_i$ associated with $F$. As was already
mentioned, in such a situation one may consider several extremal problems for poristic Steiner $4$-chains. One of the
most natural and geometrically significant target functions to be optimized is the sum of areas of circles in Steiner
chain which was studied from this point of view in a recent paper \cite{kir} (cf. also \cite{guvo}, \cite{gaz}). 
For brevity, the sum of areas of the circles in a Steiner chain will be referred to as the area of Steiner chain. In
section we determine $4$-chains having extremal area and compute those extremal values. This problem was solved
\cite{kir} in a concrete case using M\"obius transformations, but no hints have been given how to deal with the general 
case. We give a detailed solution for arbitrary Steiner $4$-chains based on the formulas for the invariants given in Section 2.

We denote by $L = \sum r_i$ the sum of poristic radii and by $S = \sum r_i^2$ the sum of poristic radii squared.
Notice that $L$ and $S$ are constant multiples of the perimeter of the quadrilateral formed by the centres of  $\Gamma_i$ and
the area of poristic $4$-chain, respectively. As above, we denote by $I_k^{(n)} = \sum b_i^k$ the $k$-th moment of
poristic curvatures. In terms of curvatures $b_i$ the first two target functions read
\begin{equation} \label{perimeter}
	L = \sum \frac{1}{b_i}.
\end{equation}
\begin{equation} \label{area}
	S = \sum \frac{1}{b_i^2}.
\end{equation}

Clearly, $L, S, I_k^{(n)}$ are symmetric functions of curvatures, which is crucial for our approach.
For further use we give the ranges of curvatures $b_i$ and their inverses $r_i$.

To study the extremal problem for area of Steiner $4$-chains we use a straightforward approach argument which was
already applied in the case of Steiner $3$-chains \cite{bisa}. The equations (\ref{moment41}), (\ref{moment42})
and Newton identities for power sums imply that the values of any symmetric function $f$ of curvatures $b_i$
are uniquely determined by the value of $b_1 \in [b_*, b^*]$. Using formulas (\ref{moment41}), (\ref{moment42}),
(\ref{moment43}) and putting $t = b_1$ we can express a symmetric target function $f$ as a function
of variable $t.$ Next, we compute the derivative $f'$ and consider its roots which yield the critical points.
It turns out that, for $S$ and $L$, the critical values of $t=b_1$ are given by explicit algebraic formulas, which 
show that the extremal values are attained exactly at the three symmetric configurations. In the case of area 
this yields the following theorem.

\begin{thm} \label{areaoptimisation}
	Extremal values of area are attained at the symmetric Steiner $4$-chains: maximum is at the axial chain $C_*^*$ and equals \\
	$A^* = [\frac{A^4 + 6A^3a + 18A^2a^2 + 6Aa^3 + a^4}{A^2 a^2 (a + A)^2}]\pi$, while the minimum value is attained at the two lateral $4$-chains and equals $A_* = [\frac{32 (3a + A)(a + 3A)}{(-a + A)^4}]\pi$
\end{thm}

{\bf Proof.} 

We take one of the curvatures, say $b_1$, as a variable denoted by $t$ and rewrite function $S$ as follows:
$$S = \frac{1}{t^2} + \frac{(b_2b_3)^2 + (b_2b_4)^2 + (b_3b_4)^2}{(b_2b_3b_4)^2}.$$
We also have 
$$b_2 + b_3 + b_4 = I_1 - t, b_2^2 + b_3^2 + b_4^2 = I_2 - t^2, b_2^3 + b_3^3 + b_4^3 = I_3 - t^3.$$

Using the well known relations between power sums and elementary symmetric functions of $b_2, b_3, b_4$ we get 

\begin{dmath}
(b_2b_3)^2 + (b_2b_4)^2 + (b_3b_4)^2 = \frac{1}{12}[3(b_2^2 + b_3^2 + b_4^2)^2 + 6(b_2 + b_3 + b_4)^2(b_2^2 + b_3^2 + b_4^2)-8(b_2 + b_3 + b_4)(b_2^3 + b_3^3 + b_4^3) - (b_2 + b_3 + b_4)^4] = \frac{1}{12}[3(I_2-t^2)^2 + 6(I_1-t)^2(I_2-t^2) - 8(I_1-t)(I_3-t^3)-(I_1-t)^4]
\end{dmath}

and

\begin{dmath}
b_2b_3b_4 = \frac{1}{6}[(b_2 + b_3 + b_4)^3 - 3(b_2 + b_3 + b_4)(b_2^2 + b_3^2 + b_4^2) + 2(b_2^3 + b_3^3 + b_4^3)]= \frac{1}{6}[(I_1-t)^3 - 3(I_1-t)(I_2-t^2) + 2(I_3-t^3),]
\end{dmath}

so that 

\begin{dmath}
S(t) = \frac{1}{t^2} + \frac{3(I_2-t^2)^2 + 6(I_1-t)^2(I_2-t^2) - 8(I_1-t)(I_3-t^3)-(I_1-t)^4}{\frac{1}{3}[(I_1-t)^3 - 3(I_1-t)(I_2-t^2) + 2(I_3-t^3)]^2}
\end{dmath}

Substituting the values of $I_1, I_2, I_3$ given by formulas (\ref{moment41}, \ref{moment42}, \ref{moment43}) we have 

\begin{equation}
	\begin{split}
	S(t) = \frac{1}{t^2}+\frac {-16{t}^{4}+64(A+a){t}^{3}-64(A+a)^{2}{t}^{2}+8({A}^{3}-{A}^{2}a-A{a}^{2}+{a}^{3})t}{(-t+A+a)^{2}((A-a)^{2}-4At-4at+4{t}^{2})^{2}}+ \\ {\frac{(A+3a)^{2}(3A+a)^{2}}{(-t+A+a)^{2}((A-a)^{2}-4At-4at+4{t}^{2})^{2}}}
	\end{split}
\end{equation}

After differentiating the above expression by $t$ and simplifying we finally get

\begin{dmath}
	\frac{dS}{dt} = -2\frac{P_1(t)P_2(t)P_4(t)}{[t(a+A-t)[(a-A)^2-4(a+A)t + 4t^2]]^3}
\end{dmath}
where
\begin{dmath*}
P_1(t) = a+A-2t,
\end{dmath*}
\begin{dmath*}
	P_2(t) = (a-A)^2 - 8(a+A)t + 8t^2,
\end{dmath*}
\begin{dmath*}
	P_4(t) = (a-A)^4(a+A)^2 - (a-A)^2(a+A)(5a^2 + 6aA + 5A^2)t + (5a^2 + 6aA + 5A^2)^2t^2-8(a+A)(5a^2 + 6aA + 5A^2)t^3 + 4(5a^2 + 6aA + 5A^2)t^4.
\end{dmath*}

The critical points of $S$ are given by the roots of those three polynomials, which are given by explicit formulas: \\
the root of $P_1$ is
$$t = \frac{a+A}{2},$$
the roots of $P_2$ are 
$$t = \frac{a+A}{2} \pm \frac{1}{2\sqrt{2}} \sqrt{a^2 + 6aA + A^2,}$$
and the roots of $P_4$ are
\begin{dmath*}
	t=\frac{A + a}{2} \pm \frac{1}{4}\sqrt{2}\sqrt{A^2 + 6A a +a^2 \pm i (A- a)^2 \sqrt{\frac{w}{5A^2 + 6A a + 5a^2}}},
\end{dmath*}

where $w=-55A^4 + 196A^3 a + 266A^2 a^2 + 196A a^3 + 55a^4$ \\

Comparing the above roots with the explicit formulas (\ref{axialcurvatures}), (\ref{offaxialcurvatures}) for the curvatures in the axial and off-axial symmetric configurations, we conclude that extremal values of area are attained exactly at those three symmetric $4$-chains. It is then easy to check that the axial $4$-chain has the maximal area, while the minimum of area is attained at the two off-axial configurations, which completes the proof of Theorem \ref{areaoptimisation}. \\

We now give examples using radii from Kiradjiev's article.\cite{kir}

\begin{ex}
	For given radii $(0.08318,0.15713,1.41421,0.15713)$, where $A$ and $a$ are equal to $-0.63585$ and $13.36461$ respectively,
	$S^*=6.461016504$.
\end{ex}	
\begin{ex}
	For given radii $(0.08318,0.15713,1.41421,0.15713)$, where $A$ and $a$ are equal to $-0.63585$ and $13.36461$ respectively,
	$S_*=1.182274825$.
\end{ex}

\begin{rem}
	Theorem \ref{areaoptimisation} also yields the interval of monotonicity of function $S$, which is important 
	in some aspects of the feasibility problem for the radii of Steiner chains. 
\end{rem}

The same argument can be used to solve several other extremal problems for Steiner $4$-chains.
As an illustration we find the extremal values of function $L$ equal to the perimeter of 
the polygon formed by the centres of circles in poristic $4$-chains.

\begin{thm} \label{perimeteroptimisation}
	Extremal values of function $L$ are attained at the symmetric Steiner $4$-chains: maximum is at the axial chain 
	$C_*^*$ and equals $L^* = -{\frac {{A}^{2}-2\,Aa+{a}^{2}}{Aa \left( a+A \right) }}$, while the minimum value is attained at the two lateral $4$-chains and equals 
	$L_* = \frac{16(A + a)}{A^2 - 2Aa + a^2}$
\end{thm}

{\bf Proof.}

We take one of the curvatures, say $b_1$, as a variable denoted by $t$ again and rewrite function $L$ as follows:
$$L = \frac{1}{t} + \frac{b_2b_3 + b_2b_4 + b_3b_4}{b_2b_3b_4}.$$
We also have 
$$b_2 + b_3 + b_4 = I_1 - t, b_2^2 + b_3^2 + b_4^2 = I_2 - t^2, b_2^3 + b_3^3 + b_4^3 = I_3 - t^3.$$
Using the well known relations between power sums and elementary symmetric functions of $b_2, b_3, b_4$ and inserting curvatures of $R$ and $r$ we get 
\begin{equation*}
	\begin{split}
		b_2b_3 + b_2b_4 + b_3b_4 =\frac{5}{4}{A}^{2}+\frac{3}{2}Aa-2At+\frac{5}{4}{a}^{2}-2at+{t}^{2}
	\end{split}
\end{equation*}
and
\begin{dmath*}
b_2b_3b_4 =\frac{1}{4}{A}^{3}+\frac{1}{4}{a}^{3}-\frac{1}{4}{A}^{2}a-\frac{5}{4}{A}^{2}t-\frac{1}{4}A{a}^{2}+2A{t}^{2}-\frac{5}{4}{a}^{2}t+2a{t}^{2}-\frac{3}{2}Aat-{t}^{3}
\end{dmath*}
so that 
$$L(t)=\frac{(A + a)(A - a)^2}
{t(-t + A + a)\left(A^2 - 2Aa - 4At + a^2 - 4a t + 4t^2\right)}$$
after differentiating by $t$ we get	
$$\frac{dL}{dt} = - \frac{(A + a)(A - a)^2 \left(A^3 - A^2 a - 10A^2 t - A a^2 - 12A a t + 24A t^2 + a^3 - 10a^2 t + 24a t^2 - 16t^3\right)}{t^2 (-t + A + a)^2 \left( A^2 - 2Aa - 4At + a^2 - 4a t + 4t^2 \right)^2}$$ 
We now put this equal to zero and solve for $t$. We get the following equation

$$-16t^3 + (24A + 24a)t^2 + (-10A^2 - 12A a - 10a^2)t + A^3 - A^2 a - A a^2 + a^3 = 0$$ 
and the critical points are the roots: \\

$t = \frac{a+A}{2},$ $t = \frac{a+A}{2} \pm \frac{1}{2\sqrt{2}} \sqrt{a^2 + 6aA + A^2,}$ \\

The proof of Theorem \ref{perimeteroptimisation} is complete.	\\
    \begin{ex}
	For given radii $(0.08318,0.15713,1.41421,0.15713)$, where $A$ and $a$ are equal to $-0.63585$ and $13.36461$ respectively,
	$L^*=1.812122475$.
    \end{ex}	
   \begin{ex}
	For given radii $(0.08318,0.15713,1.41421,0.15713)$, where $A$ and $a$ are equal to $-0.63585$ and $13.36461$ respectively,
	$L_*=1.039014169$.
   \end{ex}
\begin{rem}
	It is remarkable that the critical chains are again the symmetric chains. The fact that those chains are critical is not surprising because 
	it follows from their symmetry. So the main point of the theorem is that there are no other critical points. Analogous results can also be 
	proved for other symmetric target functions and bigger values of $n$. 
	Below we say more on the $n=6$ case. 
\end{rem}

\section*{\textbf{Concluding remarks}}

First of all, it is natural to search for generalizations of our results to Steiner $n$-chains for arbitrary $n$. For small $n$ this seems feasible in the same way as above using (\ref{yiuquadratics}) and the invariance of the first $n-1$ moments of curvatures \cite{scta}.

For $n=6$, this is especially simple since one can find the invariant first five moments
in a similar way by finding the curvatures in the axial $6$-chain by merely using \ref{eqYiu}. 
twice - for the neighbours of the maximal circle and neighbors of the minimal circle. We give 
now the main details of this computation for $n=6$.

For $n=6$, one has $q=\frac{1}{3}$ and so the coefficients of \ref{eqYiu} take the form:

\begin{equation} \label{yiu6}
	\alpha = \frac{16}{9}R^2r^2u^2, 
	\beta = -\frac{8}{3}Rru[(\frac{2}{3}Rr + (R-r)s)],
	\gamma = [\frac{4}{3}Rr - (R-r)s]^2 + 4Rrs^2.
\end{equation}

For the biggest circle in $C_*^*$ we as always have 
$s = r^* = \frac{R-r+d}{2}, v = \frac{2}{R+d-r}$ and for the smallest circle in $C_*^*$ , $s = r^* = \frac{R-r-d}{2}$ , so the sum of curvatures of their adjacent circles will be calculated as in case of $n=4$. Summing their powers one obtains explicit formulas for the moments of first five invariant curvatures. \\
The resulting formulas are rather lengthy and will be published elsewhere. \\

It is also possible to connect Theorem \ref{areaoptimisation} with similar results for Poncelet porism as follows.
It is known that the centres of poristic Steiner $n$-chains form a poristic family of polygons interscribed between 
the socle $\gamma$ and certain ellipse (see, e.g., \cite{ped}). Thus the solution given above could also
be obtained from an analogous result for Poncelet porism which follows from the formulas given in \cite{rad}.
Moreover, our results can be related to the extremal problem for an analog of Steiner chains inscribed in an arbelos
considered in \cite{gaz}.

Finally, \cite{rad} contains Fuss relations for $n$-pointed star polygons in Poncelet porism
with $n=5, 6, 7$. It is interesting to check if some of those Fuss relations are compatible with Pedoe relations \cite{ped}, which could give another relation to Poncelet porism.
Some of such bi-poristic situations have been found by A.Diakvnishvili and will be presented in her next publication.

	\bigskip
	
	Authors' e-mails and addresses:
	
	gunabibi@gmail.com
	
	ana.diakvnishvili.1@iliauni.edu.ge
	
	Ilia State University
	
	3/5, K.Cholokashvili Avenue, Tbilisi
	
\end{document}